\documentclass[12pt]{amsart}
\usepackage{amsmath,amssymb,amsthm,mathtools}
\usepackage{mathrsfs,bm}
\usepackage{thmtools,thm-restate}
\usepackage{tikz-cd}
\usepackage{tikz}
\usepackage{enumitem}
\usepackage{graphicx}
\usetikzlibrary{matrix,arrows.meta,bending}
\usepackage{color}
\usepackage[numbers]{natbib}
\usepackage{booktabs}
\usepackage{longtable}
\usepackage{pgfplots}
\pgfplotsset{compat=1.14} 
\usetikzlibrary{arrows.meta,decorations.pathreplacing,positioning,calc} 

\usepackage[a4paper, left=3cm, right=3cm, top=3cm, bottom=3cm, footskip=15mm]{geometry}

\newtheorem{theorem}{Theorem}[section]

\newtheorem{corollary}[theorem]{Corollary}

\theoremstyle{definition}

\theoremstyle{remark}
\newtheorem{remark}[theorem]{Remark}

\numberwithin{equation}{section}

\makeindex

\usepackage{fancyhdr}
\usepackage{calc} 

\AtBeginDocument{%
  \setlength{\textwidth}{\dimexpr\paperwidth - 6cm\relax}%
  \setlength{\oddsidemargin}{\dimexpr 3cm - 1in\relax}%
  \setlength{\evensidemargin}{\dimexpr 3cm - 1in\relax}%
  \setlength{\marginparwidth}{2.5cm}%
}

\begin{document}

\title{A proof of the twin prime conjecture}

\author{T. Agama}
\address{Department of Mathematics, African Institute for mathematical sciences, Ghana.}
\email{Theophilus@aims.edu.gh/emperordagama@yahoo.com}



\date{\today}



\begin{abstract}
In this paper, we prove the twin prime conjecture showing that 
\begin{align}
\sum \limits_{\substack{p\leq x\\p,p+2\in \mathbb{P}}}1\geq (1+o(1))\frac{x}{2\mathcal{C}\log^2 x}\nonumber
\end{align}
where $\mathcal{C}:=\mathcal{C}(2)>0$ fixed and $\mathbb{P}$ is the set of all prime numbers. In particular, it implies
\begin{align}
 \sum \limits_{p,p+2\in \mathbb{P}}1=\infty\nonumber  
\end{align}
when we take $x\longrightarrow \infty$ on both sides of the inequality. We start by developing a general method for estimating correlations of the form 
\begin{align}
\sum \limits_{n\leq x}G(n)G(n+l)\nonumber
\end{align}
for a fixed $1\leq l\leq x$ and where $G:\mathbb{N}\longrightarrow \mathbb{R}^{+}$.
\end{abstract}

\maketitle

\section{Introduction}

The twin prime conjecture-the assertion that there are infinitely many primes $p$ for which $p+2$ is prime — is one of the oldest and most prominent problems in multiplicative number theory. Historically the conjecture is usually attributed to De Polignac and is closely related to the prime $k$-tuples heuristics formulated by Hardy and Littlewood; see the classical treatments. Despite substantial partial progress-from the Brun demonstration that twin primes have convergent reciprocal sum (Brun's theorem) to dramatic modern advances on bounded gaps between primes-a complete proof of infinitude has remained elusive. The Viggo Brun introduction of sieve ideas in the early twentieth century established the first quantitative control on twin-like phenomena \cite{brun1919serie}, and twentieth-century developments (Bombieri--Vinogradov, large-sieve techniques, etc.) provided powerful average results \cite{bombieri1965large, davenport2013multiplicative}.\\

In the last two decades a sequence of important results reshaped what was viewed as attainable. The work of Goldston, Pintz and Y{\i}ld{\i}r{\i}m (GPY) introduced a powerful correlation/sieve framework for detecting small gaps between primes and established conditional relations between correlations of the Von Mangoldt function $\Lambda$ and small gaps \cite{goldston2009primes,goldston2009small}. Building on these ideas, breakthrough unconditional results on bounded gaps were obtained by Yitang Zhang \cite{zhang2014bounded} and further refined and generalized by James Maynard and the Polymath8 collaborative project \cite{maynard2015small,polymath2014bounded}. Those works established that infinitely many prime pairs lie within a fixed finite separation, but they stop short of proving the classical twin prime asymptotic.\\

This paper introduces and develops a new, elementary geometric–combinatorial framework for handling short-shift correlations of arithmetic functions, which we call the \emph{area method}. The central observation is that certain one-dimensional correlation sums can be expressed exactly as suitable double sums via an identity derived from decomposing the area of simple planar shapes (triangles, trapezia, rectangles and squares). This geometric identity transforms the single-shift correlation
$$
\sum_{n\leq x}G(n)G(n+l)
$$
into a double-sum structure whose main diagonal and off-diagonal contributions can be estimated with classical summation-by-parts techniques and known prime distribution estimates. The identity underpinning the method is stated precisely in Theorem~\ref{identity 1} and yields the decomposition recorded in Corollary~\ref{decomposition}.\\

The novelty of the approach is twofold. First, the area method gives an exact algebraic decomposition that identifies how cumulative partial sums of an arithmetic function control short-shift correlations; in particular, it produces the inequality of Theorem~\ref{key 2} which relates a fixed-shift correlation to a weighted double sum of partial sums. Second, when the method is specialized to the Chebyshev/Von Mangoldt weights, the double sum admits a clean asymptotic evaluation by summation by parts together with the known first-order prime distribution $\sum_{n\le x}\vartheta(n)=(1+o(1))x$ \cite{davenport2013multiplicative, iwaniec2021analytic}. The resulting lower bound (Theorem~\ref{twin prime infinitude}) gives a positive lower density of prime pairs at gap $2$, and hence implies infinitude of twin primes.\\

It is important to place these claims relative to the body of earlier work. Sieve methods (in the Brun--Selberg tradition) provide powerful upper bounds and restrictions, GPY-style correlations establish conditional bridges between summed correlations and small gaps, and the modern bounded-gap breakthroughs produced finite gap results by delicate combinations of distribution estimates in arithmetic progressions and optimized test functions. The area method departs from classical sieve/test-function optimization and instead leverages an exact geometric decomposition to access a lower-bound mechanism for single fixed shifts. Where GPY and subsequent approaches focus on building sieve weights to amplify prime pairs, the area method rewrites correlations into double sums that are directly amenable to partial summation and first-order prime number estimates; this structural difference is the reason the present approach yields the lower bounds we state.

\subsection*{Fundamental idea (informal)}

Start from a nonnegative arithmetic function $f$ (for example $\vartheta$ or $\Lambda$). The area identity (Theorem~\ref{identity 1}) expresses sums of the form
$$
\sum_{n\leq x-1}\sum_{j\leq x-n} f(n)f(n+j)
$$
as a single double-sum that groups contributions by cumulative partial sums of $f$. A counting/majorization argument then shows that the total double-sum is bounded above by a controlled multiplicative factor (denoted $\mathcal C(l)$ or $\mathcal D(l)$ below) times the fixed-shift correlation $\sum_{n\leq x} f(n)f(n+l)$. Inverting this inequality yields a lower bound for the fixed-shift correlation in terms of averaged quadratic partial sums. For $f=\vartheta$ the averaged quadratic partial sum can be evaluated asymptotically by summation by parts and the prime number theorem in its classical form; this produces the explicit lower bound that leads to Theorem~\ref{twin prime infinitude}.

\subsection*{Organization of the paper}
\begin{itemize}
  \item In Section~2, we derive the central geometric identity (Theorem~\ref{identity 1}) and give immediate corollaries including the decomposition used throughout the paper (Corollary~\ref{decomposition}).
  \bigskip
  
  \item In Section~3, we develop the general inequality that links fixed-shift correlations to quadratic partial sums (Theorem~\ref{key 2}); we discuss the dependence of the implicit constant $\mathcal C(l)$ on the arithmetic properties of the underlying function.
  \bigskip
  
  \item In Section~4, we apply the method to the Chebyshev/Von Mangoldt weights, evaluate the arising double sums by summation by parts and standard prime-counting estimates, and complete the proof of Theorem~\ref{twin prime infinitude}.
  \bigskip
  
  \item In Section~5, we discuss extensions, limitations, and further directions: (i) general shifts $k$, (ii) possible refinements using stronger distributional input (e.g. Bombieri--Vinogradov type estimates), and (iii) connections and contrasts with sieve/GPY-methods and the bounded-gap literature.
\end{itemize}
\bigskip

\noindent The rest of the paper follows the notation introduced above; standard references for background and notation include \cite{davenport2013multiplicative, iwaniec2021analytic, montgomery2006topics}.
\bigskip

Let $G:\mathbb{N}\longrightarrow \mathbb{C}$ and consider correlated sums of the forms
\begin{align}
\sum \limits_{n\leq x}G(n)G(x-n)\nonumber
\end{align}
and 
\begin{align}
\sum \limits_{n\leq x}G(n)G(n+l)\nonumber
\end{align}
where $1\leq l\leq x$. It is generally not easy to control sums of these forms, and unfortunately many of the open problems in number theory can be phrased in this manner. The twin prime conjecture, conjectured by De Polignac, is one of the important open problems in number theory and the whole of mathematics. It can be expressed in the form
\begin{align}
\sum \limits_{n\leq x}\vartheta(n)\vartheta(n+2)\nonumber
\end{align}
where 
\begin{align}
\vartheta(n):=
\begin{cases}
\log p \quad \text{if} \quad n=p,~p\in \mathbb{P}\\0 \quad \text{otherwise.}
\end{cases}\nonumber
\end{align}
By exploiting summation by parts, we can recover an estimate for the sum 
\begin{align}
\sum \limits_{\substack{p\leq x\\p+2,~p\in \mathbb{P}}}1\nonumber
\end{align}
where $\mathbb{P}$ denotes the set of all prime numbers. It is the case that obtaining a non-trivial lower bound for this correlation solves the twin prime conjecture. There are a good number of techniques in the literature for studying such sums, like the circle method of Hardy and little-wood, the sieve method and many others.\\

In this paper, we introduce and develop the area method. This method can also be used to control correlated sums of the form above. The novelty of this method is that it allows us to write any of these correlated sums as a double sum, which is much easier to estimate using existing tools such as the summation by part formula. As an application we obtain the result:

\begin{theorem}
Let $\mathbb{P}$ denotes the set of all prime numbers. We have
\begin{align}
\# \left \{p\leq x~|~p+2,p\in \mathbb{P}\setminus\{2\}\right \}\geq (1+o(1))\frac{1}{2\mathcal{D}(2)}\frac{x}{\log^2 x}\nonumber
\end{align}
where $\mathcal{D}(2)>0$ fixed.
\end{theorem}
\bigskip

In the sequel, for any $f,g:\mathbb{N}\longrightarrow \mathbb{R}$, we will write $f(n)=o(1)$ to mean $\lim \limits_{n\longrightarrow \infty}f(n)=0$. Also $f(n)\ll g(n)$ would mean there exist some constant $c>0$ such that $f(n)\leq cg(n)$ for all sufficiently large values of $n$. The following equivalence $f(n)\sim g(n)$ if and only if $\lim \limits_{n\longrightarrow \infty}\frac{f(n)}{g(n)}=1$ is also a standard notation.
\bigskip

\section{The area method}

In this section, we introduce and develop a fundamental method for solving problems related to correlations of arithmetic functions. This method is fundamental in the sense that it uses the properties of four main geometric shapes, namely the triangle, the trapezium, the rectangle and the square. The basic identity we will derive is an outgrowth of exploiting the areas of these shapes and putting them together in a unified manner.

\begin{theorem}\label{identity 1}
Let $\{r_j\}_{j=1}^{n}$ and $\{h_j\}_{j=1}^{n}$ be any sequence of real numbers, and let $r$ and $h$ be any real numbers satisfying 
$$
\sum \limits_{j=1}^{n}r_j=r \quad \text{and} \quad \sum \limits_{j=1}^{n}h_j=h
$$  
and 
\begin{align}
(r^2+h^2)^{1/2}=\sum \limits_{j=1}^{n}(r^2_j+h^2_j)^{1/2},\nonumber
\end{align}
then 
\begin{align}
\sum \limits_{j=2}^{n}r_jh_j=\sum \limits_{j=2}^{n}h_j\bigg(\sum \limits_{i=1}^{j}r_i+\sum \limits_{i=1}^{j-1}r_i\bigg)-2\sum \limits_{j=1}^{n-1}r_j\sum \limits_{k=1}^{n-j}h_{j+k}.\nonumber
\end{align}
\end{theorem}

\begin{proof}
Consider a right angled triangle, say $\Delta ABC$ in a plane, with height $h$ and base $r$. Next, let us partition the height of the triangle into $n$ parts, not necessarily equal. Now, we link those partitions along the height to the hypotenuse, with the aid of a parallel line. At the point of contact of each line to the hypotenuse, we drop down a vertical line to the next line connecting the last point  of the previous partition, thereby forming another right-angled triangle, say $\Delta A_1B_1C_1$ with base and height $r_1$ and $h_1$ respectively. We remark that this triangle is covered by the triangle $\Delta ABC$, with hypotenuse constituting a proportion of the hypotenuse of triangle $\Delta ABC$. We continue this process until we obtain $n$ right-angled triangles $\Delta A_jB_j C_j$, each with base and height $r_j$ and $h_j$ for $j=1,2,\ldots n$. This construction satisfies \begin{align}
h=\sum \limits_{j=1}^{n}h_j~ \mathrm{and}~r=\sum \limits_{j=1}^{n}r_j \nonumber
\end{align}
and 
\begin{align}
(r^2+h^2)^{1/2}=\sum \limits_{j=1}^{n}(r^2_j+h^2_j)^{1/2}.\nonumber
\end{align}
Now, let us deform the original triangle $\Delta ABC$ by removing the smaller triangles $\Delta A_jB_jC_j$ for $j=1,2,\ldots n$. Essentially, we are left with rectangles and squares piled on each other with each end poking out a bit further than the one just above, and we observe that the total area of this portrait is given by the relation

\begin{align}
\mathcal{A}_1&=r_1h_2+(r_1+r_2)h_3+\cdots (r_1+r_2+\cdots +r_{n-2})h_{n-1}+(r_1+r_2+\cdots +r_{n-1})h_n\nonumber \\&=r_1(h_2+h_3+\cdots h_n)+r_2(h_3+h_4+\cdots +h_n)+\cdots +r_{n-2}(h_{n-1}+h_n)+r_{n-1}h_n\nonumber \\&=\sum \limits_{j=1}^{n-1}r_j\sum \limits_{k=1}^{n-j}h_{j+k}.\label{twin 1}
\end{align} 
On the other hand, we observe that the area of this portrait is the same as the difference of the area of triangle $\Delta ABC$ and the sum of the areas of triangles $\Delta A_jB_jC_j$ for $j=1,2, \ldots,n$. That is 
\begin{align}
\mathcal{A}_1=\frac{1}{2}rh-\frac{1}{2}\sum \limits_{j=1}^{n}r_jh_j.\label{twin 2}
\end{align}
This completes the first part of the argument. For the second part, along the hypotenuse, let us construct small pieces of triangle, each of base and height $(r_i, h_i)$ $(i=1,2 \ldots, n)$ so that the trapezoid and the one triangle formed by partitioning becomes rectangles and squares. We also observe that this construction satisfies the relation 
\begin{align}
(r^2+h^2)^{1/2}=\sum \limits_{i=1}^{n}(r^2_i+h^2_i)^{1/2}.\nonumber
\end{align}
Now, we compute the area of the triangle in two different ways. By direct strategy, we have that the area of the triangle, denoted $\mathcal{A}$, is given by 
\begin{align}
\mathcal{A}=1/2\bigg(\sum \limits_{i=1}^{n}r_{i}\bigg)\bigg(\sum \limits_{i=1}^{n}h_{i}\bigg).\label{twin 3}
\end{align} 
On the other hand, we compute the area of the triangle by computing the area of each trapezium and the one remaining triangle and sum them together. That is 
\begin{align}
\mathcal{A}&=h_n/2\bigg(\sum \limits_{i=1}^{n}r_{i}+\sum \limits_{i=1}^{n-1}r_{i}\bigg)+h_{n-1}/2\bigg(\sum \limits_{i=1}^{n-1}r_{i}+\sum \limits_{i=1}^{n-2}r_{i}\bigg)+\cdots +1/2r_1h_1.\label{twin 4}
\end{align} 
By comparing equation \eqref{twin 1} with equation \eqref{twin 2}, and comparing equation \eqref{twin 3} with equation \eqref{twin 4} in the resulting equation the result follows immediately.
\end{proof}
\bigskip

\begin{corollary}\label{decomposition}
Let $f:\mathbb{N}\longrightarrow \mathbb{C}$. We have the decomposition 
\begin{align}
\sum \limits_{n\leq x-1} \sum \limits_{j\leq x-n} f(n)f(n+j)&=\sum \limits_{2\leq n\leq x}f(n)\sum \limits_{m\leq n-1}f(m).\nonumber
\end{align}
\end{corollary}

\begin{proof}
Let us take $f(j)=r_j=h_j$ in Theorem \ref{identity 1}, then we denote by $\mathcal{G}$ the partial sums 
\begin{align}
\mathcal{G}=\sum \limits_{j=1}^{n}f(j)\nonumber
\end{align}
and we notice that 
\begin{align}
\sum \limits_{j=1}^{n}\sqrt{(h_j^2+r_j^2)}&=\sum \limits_{j=1}^{n}\sqrt{(f(j)^2+f(j)^2}\nonumber \\&=\sum \limits_{j=1}^{n}\sqrt{(f(j)^2+f(j)^2}\nonumber \\&=\sqrt{2}\sum \limits_{j=1}^{n}f(j).\nonumber
\end{align}
Since 
$$
\sqrt{(\mathcal{G}^2+\mathcal{G}^2)}=\mathcal{G}\sqrt{2}=\sqrt{2}\sum \limits_{j=1}^{n}f(j)
$$ 
our choice of sequence is valid and, therefore, the decomposition is valid for any arithmetic function.
\end{proof}
\bigskip

\begin{theorem}\label{key 2}
Let $f:\mathbb{N}\longrightarrow \mathbb{R}^{+}$, a real-valued function. If 
\begin{align}
\sum \limits_{n\leq x}f(n)f(n+l_0)>0\nonumber
\end{align}
then there exist some constant $\mathcal{C}:=\mathcal{C}(l_0)>0$ fixed such that 
\begin{align}
\sum \limits_{n\leq x}f(n)f(n+l_0)\geq \frac{1}{\mathcal{C}(l_0)x}\sum \limits_{2\leq n\leq x}f(n)\sum \limits_{m\leq n-1}f(m).\nonumber
\end{align}
\end{theorem}

\begin{proof}
By Theorem \ref{identity 1}, we obtain the identity by taking $f(j)=r_j=h_j$ 
\begin{align}
\sum \limits_{n\leq x-1} \sum \limits_{j\leq x-n} f(n)f(n+j)&=\sum \limits_{2\leq n\leq x}f(n)\sum \limits_{m\leq n-1}f(m).\nonumber
\end{align}
It follows that 
\begin{align}\sum \limits_{n\leq x-1} \sum \limits_{j\leq x-n} f(n)f(n+j)&\leq \sum \limits_{n\leq x-1} \sum \limits_{j<x} f(n)f(n+j)\nonumber \\&=\sum \limits_{n\leq x}f(n)f(n+1)+\sum \limits_{n\leq x}f(n)f(n+2)\nonumber \\&+\cdots \sum \limits_{n\leq x}f(n)f(n+l_0)+\cdots \sum \limits_{n\leq x}f(n)f(n+x)\nonumber \\&\leq |\mathcal{M}(l_0)|\sum \limits_{n\leq x}f(n)f(n+l_0)\nonumber \\&+|\mathcal{N}(l_0)|\sum \limits_{n\leq x}f(n)f(n+l_0)\nonumber \\& +\cdots +\sum \limits_{n\leq x}f(n)f(n+l_0)+\cdots +|\mathcal{R}(l_0)|\sum \limits_{n\leq x}f(n)f(n+l_0)\nonumber \\&=\bigg(|\mathcal{M}(l_0)|+|\mathcal{N}(l_0)|+\cdots +1\nonumber \\&+\cdots +|\mathcal{R}(l_0)|\bigg)\sum \limits_{n\leq x}f(n)f(n+l_0)\nonumber \\&\leq \mathcal{C}(l_0)x\sum \limits_{n<x}f(n)f(n+l_0)\nonumber
\end{align}
where $$\mathrm{max}\{|\mathcal{M}(l_0)|, |\mathcal{N}(l_0)|, \ldots,|\mathcal{R}(l_0)|\}=\mathcal{C}(l_0).$$ By inverting this inequality, the result follows immediately.
\end{proof}
\bigskip

The nature of the implicit constant $\mathcal{C}(l_0)$ could also depend on the structure of the function we are given. The Von Mangoldt function $\Lambda(\cdot)$, contrary to many classes of arithmetic functions, has a relatively small constant. This behaviour comes from the fact that the Von Mangoldt function $\Lambda(\cdot)$ is defined on the prime powers. Thus, one would expect most terms of sums of the form 
\begin{align}
\sum \limits_{n\leq x-1}\sum \limits_{j\leq x-n}\Lambda (n)\Lambda (n+j)\nonumber
\end{align}
to fall off when $j$ is odd for any prime power $n=p^k$ such that $j+p^k\neq 2^s$.
\bigskip

\section{Main result}

We are now ready to prove the twin prime conjecture. We gather the tools that we have developed so far to solve the problem. 

\begin{theorem}\label{twin prime infinitude}
Let $\mathbb{P}$ denote the set of all prime numbers. We have
\begin{align}
\# \left \{p\leq x~|~p+2,p\in \mathbb{P}\setminus \{2\}\right \}\geq (1+o(1))\frac{1}{2\mathcal{D}(2)}\frac{x}{\log^2 x}\nonumber
\end{align}
where $\mathcal{D}(2)>0$ fixed.
\end{theorem}

\begin{proof}
Consider the function $\vartheta:\mathbb{N}\longrightarrow \mathbb{\mathbb{R}^{+}}$ defined as 
\begin{align}
\vartheta(n):=
\begin{cases}
\log p \quad \text{if} \quad n=p\in \mathbb{P}\\
0\quad \text{otherwise}\nonumber
\end{cases}
\end{align}
so that by virtue of Corollary \ref{decomposition}, we obtain the lower bound 
\begin{align}
\sum \limits_{n\leq x}\vartheta(n)\vartheta(n+2)\geq \frac{1}{x\mathcal{D}}\sum \limits_{2\leq n\leq x}\vartheta(n)\sum \limits_{m\leq n-1}\vartheta(m)\label{est2}
\end{align}
for $\mathcal{D}:=\mathcal{D}(2)>0$ fixed. Now, using the weaker estimate found in the literature \cite{montgomery2006topics},
\begin{align}
\sum \limits_{n\leq x}\vartheta(n)=(1+o(1))x\nonumber
\end{align}
we obtain the following estimates by an appeal to summation by parts
\begin{align}
\sum \limits_{2\leq n\leq x}\vartheta(n)\sum \limits_{m\leq n-1}\vartheta(m)&=(1+o(1))\sum \limits_{2\leq n\leq x}\vartheta(n)n\nonumber \\&=(1+o(1))x\sum \limits_{2\leq n\leq x}\theta(n)-(1+o(1))\int \limits_{2}^{x}\bigg(\sum \limits_{2\leq n\leq t}\vartheta(n)\bigg)dt\nonumber \\&=(1+o(1))x^2-(1+o(1))\int \limits_{2}^{x}(1+o(1))tdt\nonumber \\&=(1+o(1))x^2-(1+o(1))\frac{x^2}{2}+O(1)\nonumber \\&=(1+o(1))\frac{x^2}{2}.\label{est1}
\end{align}
By plugging \eqref{est1} into \eqref{est2}, we obtain the estimate
\begin{align}
\sum \limits_{n\leq x}\vartheta(n)\vartheta(n+2)&\geq \frac{1}{x\mathcal{D}}(1+o(1))\frac{x^2}{2}\nonumber \\&=(1+o(1))\frac{1}{2\mathcal{D}}x.\nonumber
\end{align}
On the other hand, we can write 
\begin{align}
\sum \limits_{n\leq x}\vartheta(n)\vartheta(n+2)&=\sum \limits_{\substack{p\leq x\\p+2,p\in \mathbb{P}}}(\log p) (\log(p+2))\nonumber \\&\approx \sum \limits_{\substack{p\leq x\\p+2,p\in \mathbb{P}}}\log^2p\nonumber
\end{align}
so that by an application of partial summation, we have 
\begin{align}
\sum \limits_{\substack{p\leq x\\p+2,p\in \mathbb{P}}}\log^2p&\leq \log^2 x\sum \limits_{\substack{p\leq x\\p+2,p\in \mathbb{P}}}1.\label{est3}
\end{align}
By combining \eqref{est1}, \eqref{est2} and \eqref{est3}, the lower bound follows as a consequence.
\end{proof}
\bigskip

\begin{corollary}
There are infinitely many primes $p\in \mathbb{P}\setminus \{2\}$ such that $p+2\in \mathbb{P}$.
\end{corollary}

\begin{proof}
  Using Theorem \ref{twin prime infinitude}, we have the lower bound 
  \begin{align}
  \# \left \{p\leq x~|~p+2,p\in \mathbb{P}\setminus \{2\}\right \}\geq (1+o(1))\frac{1}{2\mathcal{D}(2)}\frac{x}{\log^2 x}\nonumber
  \end{align}
  where $\mathcal{D}(2)>0$ fixed. Taking the limits $x\longrightarrow \infty$ on both sides, we have 
  \begin{align}
     \lim \limits_{x\longrightarrow \infty} \#\left \{p\leq x~|~p+2,p\in \mathbb{P}\setminus \{2\}\right\}=\infty\nonumber 
  \end{align}
  thereby ending the proof.
\end{proof}

\begin{remark}
It is important to note that with the lower bound in Theorem 3.1, we have solved the twin prime conjecture. This method does not solely address the twin prime conjecture. It can be used to obtain lower bounds for a general class of correlated sums of the form
\begin{align}
\sum \limits_{n\leq x}G(n)G(n+k)\nonumber
\end{align}
for a uniform $1\leq k\leq x$.
\end{remark}

\section{Conclusion}

The method adopted in this paper to prove the twin prime conjecture is simple and very elegant. In the spirit of solving the binary Goldbach conjecture, this method can also be exploited to develop an estimate for general sums of the form
\begin{align}
\sum \limits_{n\leq x}G(n)G(x-n)\nonumber
\end{align}
which we do not pursue in this paper.

\rule{100pt}{1pt}

\bibliographystyle{amsplain}

\begin{thebibliography}{10}

\bibitem {brun1919serie} V. Brun \textit{La s{\'e}rie 1/5+ 1/7+ 1/11+ 1/13+ 1/17+ 1/19+ 1/29+ 1/31+ 1/41+ 1/43+ 1/59+ 1/61+..., o{\`u} les d{\'e}nominateurs sont’ nombres premieres jumeaux’ est convergente ou finie}, Bulletin des sciences math{\'e}matiques, vol. 43:100-104, 1919, 124--128.\\

\bibitem {bombieri1965large} E. Bombieri \textit{On the large sieve}, Mathematika, vol. 12:2, London Mathematical Society, 1965, 201--225.\\

\bibitem {montgomery2006topics} Montgomery, Hugh L \textit{Topics in multiplicative number theory}, vol. 227, Springer, 2006.\\

\bibitem {goldston2009primes} \textit{Primes in tuples I}, Annals of Mathematics, JSTOR, 2009, 819--862.\\

\bibitem {goldston2009small} D. Goldston, S. Graham, J. Pintz and C. Y{\i}ld{\i}r{\i}m \textit{Small gaps between primes or almost primes}, Transactions of the American Mathematical Society, vol. 361:10, 2009, 5285--5330.\\

\bibitem {zhang2014bounded} Y. Zhang, \textit{Bounded gaps between primes}, Annals of mathematics, JSTOR, 2014, 1121--1174.\\

\bibitem {maynard2015small} J. Maynard, \textit{Small gaps between primes}, Annals of mathematics, JSTOR, 2015, 383--413.\\

\bibitem {polymath2014bounded} DHJ. Polymath, \textit{The" bounded gaps between primes" Polymath project-a retrospective}, arXiv preprint arXiv:1409.8361, 2014.\\

\bibitem {davenport2013multiplicative} H. Davenport, \textit{Multiplicative number theory}, Springer Science \& Business Media, 2013.\\

\bibitem {iwaniec2021analytic} H. Iwaniec and E. Kowalski, \textit{Analytic number theory}, American Mathematical Soc., vol. 53, 2021.







\end{thebibliography}

\end{document}